\documentclass{amsart}

\newtheorem{thm}{Theorem}[section]
\newtheorem{cor}[thm]{Corollary}
\newtheorem{lem}[thm]{Lemma}

\theoremstyle{definition}
\newtheorem{defn}[thm]{Definition}
\numberwithin{equation}{section}

\newcommand{\s}{\mathcal{S}}
\newcommand{\A}{\mathcal{A}}
\newcommand{\R}{\mathcal{R}}
\newcommand{\X}{\mathcal{X}}
\begin{document}

\title{Uniform non--amenability of free Burnside groups}

\author{ Denis V. Osin }

\address{}

\email{denis.osin@mtu-net.ru}

\thanks{This work has been supported by RFFR Grant $\#$ 99-01-00894
and by the INTAS Grant $\# $ 99-1224}



\date{Feb 25, 2001}

\dedicatory{}


\begin{abstract}
The aim of the present note is to show that free Burnside groups
of sufficiently large odd exponent are non--amenable in a certain
strong sense, more precisely, their left regular representations
are isolated from the trivial representation uniformly on finite
generating sets. This result is applied to the solution of a
strong version of the von Neumann -- Day problem concerning
amenability of groups without non--abelian free subgroups. As
another consequence, we obtain that the above--mentioned groups
are of uniform exponential growth. This answers a question of de
la Harpe \cite{dlH}.
\end{abstract}

\maketitle

\section{Introduction}

Let us consider the left regular representation of a finitely
generated group $G$ on the Hilbert space $L^2(G)$. For any finite
generating set $X$ of $G$, we define $\alpha (G, X)$ as the
maximal $\varepsilon \ge 0$ such that for any vector $v\in L^2(G)$
of norm $||v||=1$, there exists an element $x\in X$ satisfying the
inequality $$ || xv-v||\ge \varepsilon . $$  It is easy to check
that the existence of a finite generating set $X$ of $G$ such that
$\alpha (G,X)>0$, implies the inequality $\alpha (G,Y)>0$ for any
other finite generating set $Y\subseteq G$. Thus it is natural to
consider the quantity $$ \alpha (G)=\inf\limits_X \alpha (G,X),$$
where the infimum is taken over all finite generating sets of $G$.

Recall that a group $G$ is called {\it amenable} if there exists a
finitely additive measure $\mu $ on the set of all subsets of $G$
which is invariant under the left action of $G$ on itself and
satisfies $\mu (G)=1$. One of the most interesting
characterizations of amenable groups was obtained by Hulaniski
\cite{H} in terms of the left regular representations. In the case
of finitely generated groups it can be formulated as follows.

\begin{thm}[Hulaniski, \cite{H}]
A finitely generated group $G$ is amenable if and only if $\alpha
(G,X)=0$ for some (and hence for any) finite generating set $X$ of
$G$.
\end{thm}

In particular, we have $\alpha (G)=0$ for any amenable group $G$.
The question whether the equality $\alpha (G)=0$ is equivalent to
the amenability of $G$ was open until recent time. The first
examples of non--amenable finitely generated groups $G$ satisfying
$\alpha (G)=0$ have been constructed in \cite{Osin-WA}.

\begin{defn}
If $\alpha (G)>0$ for a finitely generated group $G$, we say that
$G$ is {\it uniformly non--amenable.}
\end{defn}

Recall that a group is said to be elementary if it contains a
cyclic subgroup of finite index. In the paper \cite{Sh}, Shalom
proved that any residually finite non--elementary hyperbolic group
is uniformly non--amenable. In fact, this is true for every
non--amenable hyperbolic group, not necessarily residually finite
(see Example 3.4 below). Essentially all uniformly non--amenable
groups known up to now are hyperbolic and, in particular, they
contain a non--abelian free subgroup. This rises the following
natural question, which is a stronger version of the so called von
Neumann--Day problem (see \cite{Day}).

{\bf Question 1.1.} {\it Does any uniformly non--amenable group
contains a non--abelian free subgroup?}

The main goal of this note is to show that the answer is negative
and can be obtained in the same way as the solution of the
classical von Neumann--Day problem. More precisely, we show that
(non--cyclic) free Burnside groups of sufficiently large odd
exponent are uniformly non--amenable. In contrast, we note that
finitely generated torsion groups which are not amenable nor
uniformly non--amenable are constructed in \cite{Osin-WA}. Our
main result is the following.

\begin{thm}
There exists an integer $N>0$ such that any free Burnside group
$B(m,n)$ of odd exponent $n>N$ and rank $m\ge 2$ is  uniformly
non--amenable.
\end{thm}

\begin{cor}
There exists a finitely generated uniformly non--amenable group
without non--abelian free subgroups.
\end{cor}

It should be mentioned that the first examples of non--amenable
groups without free subgroups were constructed by Ol'shanskii in
\cite{Ols1}. These groups are torsion with unbounded orders of
elements. After that, in \cite{A}, Adian proved that the free
Burnside groups $B(m,n)$ are non--amenable whenever $m\ge 2$, $n$
is odd, and $n\ge 665$. We use Adian's result to prove our Theorem
1.1.

The notion of uniform non--amenability is interesting, in
particular, in relation with exponential growth rates of finitely
generated groups. Recall that the {\it growth function}
$\gamma_G^X : \mathbb N \longrightarrow \mathbb N$ of a group $G$
with respect to a finite generating set $X$ is defined by $$\gamma
_G^X(n)=card\; \{ g\in G\; :\; |g|_X\le n\} .$$ The {\it
exponential growth rate} of $G$ with respect to $X$ is the number
$$\omega (G,X) = \lim_{n \to \infty} \sqrt[n]{\gamma _G^X(n)} $$
(the above limit always exists). The quantity $$\omega
(G)=\inf\limits_X \omega (G,X) $$ is called a {\it minimal
exponential growth rate} of $G$ (the infimum is taken over all
finite generating sets of $G $). One says that the group $G$ has
{\it exponential growth} if $\omega (G,X)>1$ for some (or,
equivalently, for any) finite generating subset $X$ of $G$. If
$\omega (G)>1$, then the group $G$ is said to have {\it uniform
exponential growth}.

This notion comes from geometry; in particular, if $G$ is a
fundamental group of a compact Riemannian manifold of unit
diameter, then $\log\omega (G)$ is a lower bound for the
topological entropy of the geodesic flow of the manifold
\cite{Manni}.

The first example of a group having non-uniform exponential growth
was recently constructed by J. Wilson \cite{Wil}. On the other
hand, there are many examples of classes of groups which are known
to have uniformly exponential growth, for example:

(a) non--elementary hyperbolic groups \cite{Koubi};

(b) free products with amalgamations $G\ast _{A=B}H$ satisfying
the condition $(|G:A|-1)(|H:B|-1)\ge 2$ and HNN--extensions $G\ast
_A$ associated with a monomorphisms $\phi _1, \phi _2: A\to G$,
where $|G:\phi _1(A)| + |G:\phi _2 (A)|\ge 3$ (see \cite{BH});

(c) one--relator groups of exponential growth \cite{GrH1-rel};

(d) solvable groups of exponential growth \cite{Osin-Sol} (the
particular case of polycyclic groups was considered independently
in \cite{Alp}); more generally, any elementary amenable group of
exponential growth has uniform exponential growth \cite{Osin-EG};

(e) linear groups of exponential growth \cite{Mozes}.

In the survey \cite{dlH}, de la Harpe asked the following:

{\bf Question 1.2.} {\it Do the free Burnside groups of
sufficiently large odd exponent have uniform exponential growth?}

Here we answer the question positively. In \cite{Sh}, Shalom
observed that every uniformly non--amenable finitely generated
group has uniform exponential growth. From this and Theorem 1.1,
we obtain immediately

\begin{cor}
If $m\ge 2$, $n$ is odd and large enough, then the group $B(m,n)$
has uniform exponential growth.
\end{cor}

It is worth to note, that in some similar sense uniform
non--amenability of free Burnside groups of large odd exponent is
claimed in the paper \cite{AetAll}. However, in \cite{AetAll}
authors uses an unproved result about free subgroups of $B(m,n)$.
This result was announced by Ivanov in \cite{IO}, but the complete
proof have never been published.

{\bf Acknowledgements.} I would like to thank Alexander Yu.
Ol'shanskii for useful comments concerning the construction of
free Burnside groups.

\section{Constructing free Burnside groups}

In this section we recall the construction of presentations of
free Burnside groups and describe shortly the main properties used
in this paper. For details we refer to the books \cite{Abook},
\cite{Obook}. Our discussion heavily depends on results from
\cite{Obook} (for otherwise our paper would be unreasonably long).
In fact, we will not obtain any new facts about Burnside groups by
analyzing the geometric structure of van Kampen diagrams; we will
only show how to combine certain lemmas from \cite{Obook} in order
to obtain the result we need. Our main goal here is to prove
Theorem 2.7.

Given an alphabet $\mathcal A$, we denote by $|W|$ the length of a
word $W$ over $\mathcal A$. For two words $U,V$ over $\mathcal A$
we write $U\equiv V$ to express letter--by--letter equality.
Finally if $\mathcal A$ is a generating set of a group $G$, we
write $U=V$ whenever two words $U$ and $V$ over $\mathcal A ^{\pm
1}$ represent the same elements of $G$; we identify the words over
$\A $ and the elements of $G$ represented by them.

Recall that $B(m,n)$, the free Burnside group of exponent $n$ and
rank $m$, is the free group in the variety defined by the low
$X^n=1$. Throughout this paper we will assume that $n$ is odd and
large enough, and $m\ge 2$; all lemmas and theorems are formulated
under these assumptions.

The group $B(m,n)$ can be defined by the presentation
\begin{equation}
B(\A, n)=\langle \A \; : \; R=1, R\in \bigcup\limits_{i=1}^\infty
\mathcal R_i\rangle ,\label{B}
\end{equation}
where $\mathcal A =\{ a_1, \ldots , a_n\} $ and the sets of
relations $\mathcal R_i$ are constructed as follows \cite{Obook}.
We put $\mathcal R_0 =\emptyset $. By induction, suppose that we
have already defined the set of relations $\mathcal R_{i-1}$,
$i\ge 1$. Denote by $G(i-1)$ the group with the presentation
$\langle \A \; : \; R=1, R\in \mathcal R_i\rangle $.

For $i\ge 1$, a word $X$ over the alphabet $\A $ is called {\it
simple in the rank} $i-1$, if it is not conjugated to a power of a
shorter word in the group $G(i-1)$ and is not conjugated to a
power of a period of rank $k\le i-1$ in the group $G(i-1)$. Let us
denote by $\mathcal X_i$ a certain maximal subset of words
satisfying the following conditions.

1) $\mathcal X_i$ consists of words of length $i$ which are simple
in the rank $i-1$.

2) If $A, B\in \mathcal X_i$ and $A\not\equiv B$, then $A$ is not
conjugated to $B$ or $B^{-1}$ in the group $G(i-1)$.

\noindent Each word from $\mathcal X_i$ is called a {\it period of
rank $i$.} We introduce the additional relations $\s _i=\{ A^n\; :
\; A\in \X _i\} $ and set $\R _i=\R _{i-1} \cup \s _i$.

Now we are going to list some results needed for the sequel.

\begin{lem}[\cite{Obook}, Theorem 19.4]\label{1}
Any element of $B(m,n)$ is conjugated to a period of a certain
rank.
\end{lem}

\begin{lem}[\cite{Obook}, Theorem 19.5]\label{2}
The centralizer of every nontrivial element of $B(m,n)$ is a
cyclic subgroup of order $n$.
\end{lem}

From Lemma \ref{2} we have immediately
\begin{cor}\label{3}
Suppose that $a, b$ are two elements of $B(m,n)$ such that
$[a,b]\ne 1$. Then $[a^i, b]\ne 1$ whenever $a^i\ne 1$.
\end{cor}

\begin{proof}
Suppose that $[a^i,b]=1$ and $a^i\ne 1$. By Lemma \ref{2}, the
centralizer of $a^i$ is cyclic. As $a, b\in C(a^i)$, we have
$a=z^{m_1}$ and $b=z^{m_2}$ for a certain $z\in B(m,n)$. This
contradicts to $[a,b]\ne 1$.
\end{proof}

In the following two lemmas the additional parameters $d$ and
$\sigma $ appear (the parameter $\sigma $ here corresponds to
$100\zeta ^{-1}$ from \cite{Obook}). Exact values of these
parameters are not important for us (in fact, $d,\sigma << n$, but
the only inequality we need here is $\sigma <\frac{n}{4}$). It is
sufficient to know that {\it there exists parameters $d$ and
$\sigma $ such that the following two lemmas hold}.

\begin{lem}\label{4}
Let $C$ be a period of a certain rank, $V\equiv C^k$, where
$\sigma <k\le \frac{1}{2} n$. Suppose that an element $W$ does not
commute with $V$ and has minimal length among elements of the
double coset $\langle C^k\rangle W \langle C^k\rangle $. We also
assume that $[C^k, W]$ is conjugated to $A^l$, where $A$ is a
period of a certain rank and $|l|\le \frac{1}{2} n$. Then we have
$|l|\le \sigma $ and the pair $[C^k, W], C^k$ is conjugated to the
pair $(A^l, B)$, where $|B|<d|A|$.
\end{lem}

{\it On the proof}. The lemma is a simplification of Lemma 25.21
from \cite{Obook}. We have to note that the presentation (\ref{B})
satisfies the condition $R$ (there are no relations of the second
type at all, see \cite[Ch. 8]{Obook} for definitions and details).

\begin{lem}\label{5}
Suppose that $l$ is an integer and $A, B\in B(m,n)$ such that

1) $A$ is a period of a certain rank;

2) $[A,B]\ne 1$ in $B(m,n)$;

3) $|l|\le \sigma $ and $|B|< d|A|$.

\noindent Then there exists integer $s$ such that the pair
$(BA^{ls}, B)$ is conjugated to a pair $(F, T)$, where $F$ is a
period of a certain rank, $[F, T]\ne 1$, and $|T|\le 3F$.
\end{lem}

{\it On the proof.} The proof of this lemma is a part of the proof
of Lemma 27.3 from \cite{Obook}. Namely in order to prove the
above assertion we have to repeat (word for word) the last 5
paragraphs of the proof of Lemma 27.3 from \cite{Obook}.

\begin{lem}\label{6}
For every odd large enough $n\in \mathbb N$, there exist $h,k\in
\mathbb N$ with the following property. Suppose that $F$ is a
period of a certain rank, $T$ is an element of $B(m,n)$ such that
$[F, T]\ne 1$, and $|T|\le 3F$. Then the elements $$
b_1=F^kTF^{k+2}\ldots TF^{k+2h-2}$$ and
$$b_1=F^{k+1}TF^{k+3}\ldots TF^{k+2h-1}$$ form the basis of the
free Burnside group of exponent $n$.
\end{lem}

\begin{proof}
Up to notation the statement of the lemma can easily be obtained
from Lemma 39.4 \cite{Obook}. We only notice that our constants
$n$, $k$, and $h$ correspond to constants $n_0$, $n$ and $h$ from
\cite{Obook} respectively.
\end{proof}

Now we are ready to prove the main result of this section.

\begin{thm}\label{th}
For any odd large enough $n$ and $m\ge 2$, there exists $M\in
\mathbb N$ having the following property. Let $a,b$ be two
non--commuting elements of $B(m,n)$. Then there exist  elements
$u,v\in \langle a,b \rangle $ such that $\{ u,v\} $ is the basis
of the free Burnside subgroup of exponent $n$ and the lengths of
elements $u,v$ with respect to $\{ a,b\} $ satisfy $$ |u|_{\{
a,b\} } < M, \;\;\;\;\;\; |v|_{\{ a,b\} }<M.$$
\end{thm}

\begin{proof}
By Lemma \ref{1}, we have $a=P^{-1}C^{k^\prime }P$, for some
element $P$ and a period $C$ of a certain rank. Without loss of
generality we can assume that $0<k^\prime < \frac{n}{2}$. Recall
that $\sigma <\frac{n}{4}$. Clearly there is a number $i\in
\mathbb N$ such that
\begin{equation}\label{T2}
i< \frac{n}{2}
\end{equation}
and $\sigma < ik^\prime <\frac{n}{2}.$ Thus we have
\begin{equation}\label{T312}
a^i=P^{-1}C^kP,
\end{equation}
where
\begin{equation}\label{T4}
\sigma <k<\frac{n}{2}.
\end{equation}
Using Corollary \ref{3}, we note that $$ [C^k, PbP^{-1}]=P[a^i,
b]P^{-1}\ne 1.$$ Denote by $W_0$ the element $PbP^{-1}$ and by $W$
the shortest element in the double coset $\langle C^k\rangle W_0
\langle C^k\rangle $. Evidently $W=C^{kk_1}W_0C^{kk_2}$ for some
$k_1, k_2$ satisfying the inequality
\begin{equation}\label{T7}
\max \{ |k_1|, |k_2|\} < \frac{n}{2}.
\end{equation}
Since $W=Pa^{ik_1}ba^{ik_2}P^{-1}$, we obtain the following
inequality by using (\ref{T7}) and (\ref{T2})
\begin{equation}\label{T8}
|P^{-1}WP|_{\{ a,b\} }\le 2i+k_1 +k_2+1< 2n.
\end{equation}
It is clear that $[C^k, W]\ne 1$. By Lemma \ref{1}, $[C^k, W]$ is
conjugated to $A^l$, where $A$ is a period of a certain rank,
i.e.,
\begin{equation}\label{***}
[C^k, W]=Q^{-1}A^lQ
\end{equation}
for some element $Q$. Set
\begin{equation}\label{**}
B=Q^{-1}C^kQ.
\end{equation}
Applying Lemma \ref{4}, we obtain $|l|\le \sigma $ and $|B|\le
d|A|$. Note that all conditions of Lemma \ref{5} are satisfied for
$A$ and $B$. Therefore, there is integer $s$ such that the pair
$(BA^{ls}, B)$ is conjugated to a pair $(F,T)$, where $F$ is a
period of a certain rank, $[F,T]\ne 1$, and $|T|\le 3|F|$.
According to Lemma \ref{6} there exists $b_1, b_2\in \langle
F,T\rangle $ such that $b_1, b_2$ freely generate the subgroup
$\langle b_1, b_2\rangle $ and $$ \max \{ |b_1|_{\{ F, T\} },
|b_2|_{\{ F, T\} }\} \le N,$$ where $N$ depends on $k$ and $h$
only (it is easy to calculate the lengths of $b_1$ and $b_2$
exactly, but this is not our goal here). Therefore, there exists
$x_1, x_2\in \langle BA^{ls}, B\rangle $ such that $x_1, x_2$
freely generate the subgroup $\langle x_1, x_2\rangle $ and $$
\max \{ |x_1|_{\{ BA^{ls}, B\} },\; |x_2|_{\{ BA^{ls}, B\} }\} \le
N.$$ Clearly we can assume $ls<\frac{n}{2}$. Passing to generators
$A$ and $B$, we obtain $$ \max \{ |x_1|_{\{ A, B\} },\; |x_2|_{\{
A, B\} }\} \le N\left( \frac{n}{2}+1\right) .$$ Therefore taking
into account (\ref{***}) and (\ref{**}), we obtain the following
estimates for the elements $y_i=Qx_iQ^{-1}$, $i=1,2$,  $$ \max \{
|y_1|_{\{ C^k, W\} },\; |y_2|_{\{ C^k, W\} }\} \le 4\max \{
|x_1|_{\{ A, B\} },\; |x_2|_{\{ A, B\} }\} \le N(2n+4).$$ Finally
we set $z_i=P^{-1}y_iP$ for $i=1,2$. Combining (\ref{T8}),
(\ref{4}), and (\ref{T312}), we obtain $$ \max \{ |z_1|_{\{ a, b\}
},\; |z_2|_{\{ a, b\} }\} \le 2n(2n+4)N.$$ As the pair $(z_1,
z_2)$ is conjugated to $(x_1, x_2)$, the elements $z_1, z_2$ form
the basis of free Burnside group of exponent $n$. To conclude the
proof it suffices to set $M=2n(2n+4)N$.
\end{proof}

{\bf Remark.} It should be noted that Ivanov announced  the
following stronger result in \cite{IO}. Given any $n$ odd and
large enough, and $m\ge 2$, there exist words $u(x,y), v(x,y)$
over the alphabet $\{ x,y\} $ such that if $a,b$ are two
non--commuting elements of $B(m,n)$, then $u(a,b), v(a,b)$
generate freely a free Burnside subgroup of $B(m,n)$.
Unfortunately the proof of this result has never been written.

\section{Sufficient conditions for uniform non--amenability}

\begin{defn} Suppose that $G$ is a group with a given finite set of
generators $X$, $Y$ a subset of $G$. The {\it depth } of $Y$ with
respect to $X$ is defined by $$depth _X(Y)=\max\limits_{y\in Y} |
y| _X. $$
\end{defn}

\begin{lem} Suppose that $G$ is a group, $H$ is a
subgroup of $G$, $X$ and $Y$ are finite generating sets of $G$ and
$H$ respectively. Then we have \begin{equation} \alpha (G,X) \ge
\frac{\alpha (H,Y)}{depth _X (Y)}.\label{aGH} \end{equation}
\end{lem}

{\it Proof. } In order to prove (\ref{aGH}), we have to show that
for any vector $f\in L^2(G)$ of norm one, there exists an element
$x\in X^{\pm 1}$ such that
\begin{equation}
\| xf - f \| \ge \frac{\alpha (H,Y)}{depth _X (Y)}. \label{fxfad}
\end{equation}

As usual we denote by $Hs$ the right coset representing by $s$.
Let us fix a (unique) element $s$ in each right coset and denote
by $S$ the obtained system of representatives. Given a function
$f\in L^2(G)$, for every $s\in S$, we introduce a new function $$
f_s(g)=\left\{ \begin{array}{ll} f(g), & {\rm if } \; g\in Hs,\\
0, & {\rm if } \; g\notin Hs. \end{array} \right. $$ If $f_s\ne
0$, we set
\begin{equation}
\tilde f_s = \frac{f_s}{\| f_s\| }.\label{40}
\end{equation}
Obviously the norm of $\tilde f_s$ is equal to $1$ whenever
$f_s\ne 0$. It is also clear that $$f=\sum\limits_{s\in S} f_x.$$
Further, to each $f_s$, $s\in S$, we assign a function $ h_s\in
L^2(H)$ by the rule $$ h_s(g)=\tilde f_s(gs).$$ Note that $\|
h_s\| =\| \tilde f_s\| =1. $ Therefore, by definition of $\alpha
(H,Y)$, there is an element $y\in Y^{\pm 1}$ such that
\begin{equation}
\| yh_s- h_s\| \ge \alpha (H,Y). \label{41}
\end{equation}
On the other hand, we have
\begin{equation}
\begin{array}{rl}
\| y \tilde f_s -\tilde f_s\| & = \sqrt{\sum\limits_{g\in Hs}
\left( \tilde f_s(y^{-1}g)-\tilde f_s(g) \right) ^2} =
\sqrt{\sum\limits_{g\in H} \left( \tilde f_s(y^{-1}gs)-\tilde
f_s(gs) \right) ^2} = \\ & \\ & \sqrt{\sum\limits_{g\in Hs} \left(
h_s(y^{-1}g)- h_s(g) \right) ^2} = \| yh_s - h_s\|. \label{42}
\end{array}
\end{equation}
Combining (\ref{40}), (\ref{41}), and (\ref{42}) yields
\begin{equation}
\| y f_s - f_s \| \ge \alpha (H, Y) \| f_s\| .\label{43}
\end{equation}

Now we observe that the supports of the functions $(yf_s - f_s)$
are pairwise disjoint. Thus $f_s$, $s\in S$, are pairwise
orthogonal and we obtain
\begin{equation}
\begin{array}{rl}
\| yf -f\| = & \left\| y\sum\limits_{s\in S} f_s -
\sum\limits_{s\in S} f_s\right\| = \sqrt{\sum\limits_{s\in S} \|
yf_s - f_s \| ^2}\ge \\ &\\ &  \alpha (H,Y)\sqrt{\sum\limits_{s\in
S} \| f_s \| ^2}=\alpha (H,Y).
\end{array} \label{44}
\end{equation}

Assume that there exists a function $f\in L^2(G)$ that does not
satisfy (\ref{fxfad}), i.e.,
\begin{equation}
\| xf - f \| < \frac{\alpha (H,Y)}{depth _X (Y)} \label{45}
\end{equation}
for every $x\in X^{\pm 1}$. According to the above--mentioned
arguments, there exists an element $y\in Y$ satisfying the
inequality (\ref{44}). Suppose that $y=x_d\ldots x_2 x_1$, where
$x_d, \ldots , x_2, x_1 \in X^{\pm 1}$ and $d\le depth _X(Y)$.
Denote by $w_i$ the $i-th$ suffix of the word $x_d\ldots x_2 x_1$,
i.e., the word $x_i\ldots x_2 x_1$, and set $w_0 = 1$ for
convenience. The inequality (\ref{45}) implies $$
\begin{array}{rl} \| yf -f\| = & \left\| \sum\limits_{i=0}^{d-1}
(w_{i+1} f -w_i f)\right\| \le \sum\limits_{i=0}^{d-1} \| x_{i+1}
w_if-w_if\| = \\ & \\ & \sum\limits_{i=0}^{d-1} \| x_{i+1} f-f\| <
d \cdot \frac{\alpha (H,Y)}{depth _X(Y)}\le \alpha (H,Y).
\end{array} $$
This contradicts to (\ref{44}). The theorem is proved. $\Box $

\begin{cor}\label{cor}
Suppose that $H$ is a non--amenable group with a given set of
generators $Y$, $G$ is a finitely generated group. Assume that for
every finite generating set $X$ of $G$ there exists an embedding
$\iota : H\to G$ such that $$depth _X(\iota (Y))\le D,$$ where the
constant $D$ is independent of $X$. Then the group $G$ is
uniformly non--amenable.
\end{cor}

{\bf Example 3.4.} Any infinite hyperbolic group is uniformly
non--amenable. Indeed, by the result of Koubi \cite{Koubi}, we can
apply the Corollary \ref{cor} for $H=F_2$, the free group of rank
$2$.

\begin{proof}[Proof of Theorem 1.1.]
Recall that the group $B(2,n)$ is non--amenable for all odd
$n>665$. To prove the theorem it is sufficient to refer to Theorem
\ref{th} and Corollary \ref{cor}.
\end{proof}


\end{document}